\documentclass[12pt]{article}

\usepackage{amsmath,amssymb,enumerate,bbm}
\newcommand{\tmop}[1]{\ensuremath{\operatorname{#1}}}

\newtheorem{theorem}{Theorem}[section]

\newtheorem{lemma}[theorem]{Lemma}

\numberwithin{equation}{section}
\newenvironment{proof}{\noindent\textbf{Proof\ }}{\hspace*{\fill}$\Box$\medskip}

 \makeatletter
 
 \newcommand{\Rmnum}[1]{\expandafter\@slowromancap\romannumeral #1@}
 \makeatother

\begin{document}

\title{On $k$-simplexes in $(2k-1)$-dimensional vector spaces over finite fields}\author{Le Anh Vinh\\
Mathematics Department\\
Harvard University\\
Cambridge, MA 02138, US\\
vinh@math.harvard.edu}\maketitle

\begin{abstract}
  We show that if the cardinality of a subset of the $(2k-1)$-dimensional vector
  space over a finite field with $q$ elements is $\gg q^{2k-1-\frac{1}{2k}}$, then it contains a positive proportional of all $k$-simplexes up to congruence.
\end{abstract}

\begin{center}
Mathematics Subject Classifications: 05C15, 05C80.\\
Keywords: finite Euclidean graphs, finite non-Euclidean graphs, pseudo-random graphs.
\end{center}

\section{Introduction}
\label{sec:in}

A classical result due to Furstenberg, Katznelson and Weiss \cite{furstenberg} says
that if $E \subset \mathbbm{R}^2$ has positive upper Lebesgue density, then
for any $\delta > 0$, the $\delta$-neighborhood of $E$ contains a congruent
copy of a sufficiently large dilate of every three point configuration. For higher dimensional simplexes, Bourgain \cite{bourgain} showed that if $E \subset \mathbbm{R}^d$ has positive upper density and $\Delta$ is a $k$-simplex with $k<d$, then $E$ contains a rotated and translated image of every large dilate of $\Delta$. The cases $k=d$ and $k=d+1$ are still remain open. Akos Magyar \cite{magyar1,magyar2} studied  related problems in the integer lattice $\mathbbm{Z}^d$. He showed \cite{magyar2} that if $d > 2k+4$ and $E \subset \mathbbm{Z}^d$ has positive upper densitiy, then all large (depending on density of $E$) dilates of a $k$-simplex in $\mathbbm{Z}^d$ can be embedded in $E$. 

Hart and Iosevich \cite{hart-iosevich} made the first investigation in an analog of this question in finite field geometries. They showed that if $E \subset \mathbbm{F}_q^d$, $d \geq \binom{k+1}{2}$, such that $|E| \geq Cq^{\frac{kd}{k+1}}q^{\frac{k}{2}}$ with a sufficiently large constant $C>0$, then $E$ contains an isometric copy of every $k$-simplex. Using graph theoretic method, the author \cite{vinh-dg} showed that the same result holds for $d  \geq 2k$ and $|E| \gg q^{\frac{d-1}{2}+k}$ (cf. Theorem 1.4 in \cite{vinh-dg}). 

Note that, serious difficulties arise when the size of the simplex is sufficiently large with respect to the ambient dimension. Even in the case of triangles, the result in \cite{vinh-dg} is only non-trivial for $d \geq 4$.   
In \cite{covert}, Covert, Hart, Iosevich and Uriarte-Tuero addressed the case of triangles in two-dimensional vector spaces over finite fields. They showed that if $E$ has density $\geq \rho$, for some $\frac{C}{\sqrt{q}} \leq \rho \leq 1$ with a sufficiently large constant $C >0$, then the set of triangles determined by $E$, up to congruence, has density $\geq c\rho$. In \cite{vinh-tfs}, the author studied the remaining case, triangles in $3$-dimensional vector spaces over finite fields. Using a combination of graph theory method and Fourier analysis, the author showed that if $E \subset \mathbbm{F}_q^d$, $d\geq 3$, such that $|E| \gg q^{\frac{d+2}{2}}$, then $E$ determines almost all triangles up to congruence. The arguments in \cite{vinh-tfs} however do not work for $k\geq 5$.  

In this paper, we will study the case of $k$-simplexes in $(2k-1)$-dimensional vector spaces with $k\geq 3$. Given $E_1,\ldots,E_k \subset \mathbbm{F}_q^d$, where $\mathbbm{F}_q$ is a finite field of
$q$ elements, define
\begin{equation}
  T_k (E_1,\ldots,E_k) =\{(x_1, \ldots x_k) \in E_1 \times \ldots \times E_k \}/ \sim
\end{equation}
with the equivalence relation $\sim$ such that $(x_1,\ldots,x_k) \sim (x_1', \ldots,x_k')$
if there exists $\tau \in \mathbbm{F}_q^d$ and $O \in S O_d (\mathbbm{F}_q)$,
the set of $d$-by-$d$ orthogonal matrices over $\mathbbm{F}_q$ with
determinant $1$, such that
\begin{equation}
  (x_1',\ldots,x_k') = (O (x_1) + \tau, \ldots, O(x_k) + \tau) .
\end{equation}

The main result of this paper is the following. 

\begin{theorem}\label{main2-tes}
Let $E \subset \mathbbm{F}^{2k-1}_q$ with $k\geq 3$, and suppose that
\[|E| \gg q^{2k-1-\frac{1}{2k}}.\]
Then there exists $c>0$ such that
\[|T_{k+1}(E)| \geq c q^{\binom{k+1}{2}}.\]
\end{theorem}

In other words, we always get a positive proportion of all $k$-simplexes if $E \gg q^{k-1-\frac{1}{2k}}$ and $k \geq 3$. The rest of this short paper is organized as follows. In Section \ref{sec:star}, we establish some results about the occurrences of colored subgraphs in a pseudo-random coloring of a graph. In Section \ref{sec:graphs}, we construct our main tools to study simplexes in vector spaces over finite fields, the finite Euclidean and non-Euclidean graphs. We then prove our main result, Theorem \ref{main2-tes} in Section \ref{sec:proofs}.

\section{Subgraphs in expanders}
\label{sec:star}

We call a graph $G = (V, E)$ $(n, d, \lambda)$-graph if $G$ is a $d$-regular
graph on $n$ vertices with the absolute values of each of its eigenvalues but
the largest one is at most $\lambda$. Suppose that a graph $G$ of order $n$ is colored by $t$ colors. Let $G_i$ be
the induced subgraph of $G$ on the $i^{\tmop{th}}$ color. We call a
$t$-colored graph $G$ $(n, d, \lambda)$-r.c. (regularly colored)
graph if $G_i$ is a $(n, d, \lambda)$-regular graph for each color $i \in \{1,
\ldots, t\}$. In this section, we will study the occurrences of colored subgraphs in $(n,d,\lambda)$-r.c. graphs. 

\subsection{Colored subgraphs}
It is well-known that if $\lambda \ll d$
then an $(n, d, \lambda)$-graph behaves similarly as a random graph $G_{n, d /
n}$. Precisely, we have the following result.

\begin{theorem} (cf. Theorem 9.2.4 in \cite{alon-spencer}) \label{tool 1}
  Let $G$ be an $(n, d, \lambda)$-graph. For a vertex $v \in V$ and a subset
  $B$ of $V$ denote by $N (v)$ the set of all neighbors of $v$ in $G$, and
  let $N_B (v) = N (v) \cap B$ denote the set of all neighbors of $v$ in $B$.
  Then for every subset $B$ of $V$:
  \begin{equation}
    \sum_{v \in V} (|N_B (v) | - \frac{d}{n} |B|)^2 \leqslant
    \frac{\lambda^2}{n} |B| (n - |B|) .
  \end{equation}
\end{theorem}

The following result is an easy corollary of Theorem \ref{tool 1}.

\begin{theorem} \label{tool 2}
  (cf. Corollary 9.2.5 in \cite{alon-spencer}) Let $G$ be an $(n, d, \lambda)$-graph. For every
  set of vertices $B$ and $C$ of $G$, we have
  \begin{equation}
    |e (B, C) - \frac{d}{n} |B\|C\| \leqslant \lambda \sqrt{|B\|C|},
  \end{equation}
  where $e (B, C)$ is the number of edges in the induced bipartite subgraph of
  $G$ on $(B, C)$ (i.e. the number of ordered pair $(u, v)$ where $u \in B$,
  $v \in C$ and $u v$ is an edge of $G$). 
\end{theorem}

Let $H$ be a fixed graph of order $s$ with $r$ edges and with automorphism
group $\tmop{Aut} (H)$. It is well-known that for every constant $p$ the
random graph $G (n, p)$ contains
\begin{equation}
  (1 + o (1)) p^r (1 - p)^{(^s_2) - r} \frac{n^s}{| \tmop{Aut} (H) |}
\end{equation}
induced copies of $H$. Alon extended this result to $(n, d, \lambda)$-graph.
He proved that every large subset of the set of vertices of a $(n, d,
\lambda)$-graph contains the ``correct'' number of copies of any fixed small
subgraph (Theorem 4.10 in \cite{krivelevich-sudakov}).

\begin{theorem}\label{tool 22}
  (\cite{krivelevich-sudakov}) Let $H$ be a fixed graph with $r$ edges, $s$
  vertices and maximum degree $\Delta$, and let $G = (V, E)$ be an $(n, d,
  \lambda)$-graph, where, say, $d \leqslant 0.9 n$. Let $m < n$ satisfies $m
  \gg \lambda \left( \frac{n}{d} \right)^{\Delta}$. Then, for every subset $U
  \subset V$ of cardinality $m$, the number of (not necessrily induced) copies
  of $H$ in $U$ is
  \begin{equation}
    (1 + o (1)) \frac{m^s}{| \tmop{Aut} (H) |} \left( \frac{d}{n} \right)^r .
  \end{equation}
\end{theorem}

In \cite{vinh-dg}, we observed that Theorem \ref{tool 22} can be extended to $(n, d, \lambda)$-r.c. graph. Precisely, we showed that every large subset of the set of vertices of an $(n, d, \lambda)$-r.c. graph contains the ``correct'' number of copies of any fixed small colored graph. We present here a multiset version of this statement. 

\begin{theorem} \label{tool 3}
  Let $H$ be a fixed $t$-colored graph with $r$ edges, $s$ vertices, maximum
  degree $\Delta$ (with the vertex set is ordered), and let $G$ be a $t$-colored graph of order $n$.
  Suppose that $G$ is an $(n, d, \lambda$)-r.c graph, where, say, $d \ll n$.
  Let $E_1,\ldots,E_s \subset V$ satisfy $|E_i| \gg \lambda \left( \frac{n}{d} \right)^{\Delta}$.
  Then the number of (not
  necessrily induced) copies of $H$ in $E_1 \times \ldots \times E_s$ (one vertex in each set) is
  \begin{equation}
    (1 + o (1)) \prod_{i=1}^s |E_i| \left( \frac{d}{n} \right)^r.
  \end{equation}
\end{theorem}

The proof of this theorem is similar to the proofs of Theorem 4.10 in
\cite{krivelevich-sudakov} and Theorem 2.3 in \cite{vinh-dg}. 
Note that going from one color formulation (Theorem 4.10 in \cite{krivelevich-sudakov}) and one set formulation (Theorem 2.3 in \cite{vinh-dg}) to a multicolor and multiset formulation (Theorem \ref{tool 3}) is just a matter of inserting different letters in couple places. 

\subsection{Colored stars}

Given any $k$ colors $r_1,\ldots,r_k$, a $k$-star of type $(r_1,\ldots,r_k)$ has $k+1$ vertices, one center vertex $x_0$ and $k$ leaves $x_1,\ldots,x_k$, with the edge $(x_0,x_i)$ is colored by the color $r_i$. The following result gives us an estimate for the number of colored $k$-star in an $(n,d,\lambda)$-r.c. graph $G$ (see \cite{vinh-tfs} for an earlier version).

\begin{theorem}\label{main-tool}
  Let $G$ be an $(n, d, \lambda)$-r.c. graph. Given any $k$ colors $r_1,\ldots,r_k$ in the color set. Suppose that $E_0,E_1,\ldots,E_k \subset V(G)$ with 
  \begin{equation}\label{big} |E_0|^2\prod_{i \in I} |E_i| \gg \left(\frac{n}{d}\lambda \right)^{2|I|}
  \end{equation}
  for all $I \subset \{1,\ldots,k\}$, $|I| \geq 2$, and
  \begin{equation}\label{small}|E_0||E_i| \gg \left(\frac{n}{d}\lambda \right)^2\end{equation}
  for all $i \in \{1,\ldots,k\}$. 
  Let $e_{\{r_1,\ldots,r_k\}}(E_0;\{E_1,\ldots,E_k\})$ denote the number of $k$-star of type $(r_1,\ldots,r_k)$ in $E_0\times E_1 \times \ldots \times E_k$ (with the center in $E_0$). Then
  \begin{equation}
    e_{\{r_1,\ldots,r_k\}}(E_0;\{E_1,\ldots,E_k\}) = (1+o(1))\left(\frac{d}{n}\right)^k\prod_{i=0}^k |E_i|,
  \end{equation}
  where $k$ is fixed and $n, d, \lambda \gg 1$.
\end{theorem}

\begin{proof}
The proof proceeds by induction. The base step is the case $k=1$. Since $|E| \gg \frac{n}{d}\lambda$ and the number of $1$-stars of type $a$ in $E_0 \times E_1$ is just the number of $a$-colored edges in $E_0 \times E_1$, the statement follows immediately from Theorem \ref{tool 2} and (\ref{small}).

Assuming that the statement holds for all colored $l$-stars with $l<k$. For a vertex $v \in V$ and a color $r$, let $N_E^{r} (v)$ denote the set of all
  $r$-colored neighbors of $v$ in $E$. From Theorem
  \ref{tool 1}, we have
  \begin{equation}\label{1}
    \sum_{v \in E_0} (|N_{E_i}^{r_i} (v) | - \frac{d}{n} |E_i|)^2 \leqslant \sum_{v \in V}
    (|N_{E_i}^{r_i} (v) | - \frac{d}{n} |E_i|)^2  \leqslant  \frac{\lambda^2}{n} |E_i|
    (n - |E_i|) \leqslant \lambda^2 |E_i|. \end{equation}
    
For $k \geqslant 2$, by the Cauchy-Schwartz inequality, we have
\begin{equation}\label{2} \prod_{i = 1}^k \left( \sum_{j = 1}^n a_{i.j}^2 \right) \geqslant \left(
   \sum_{j = 1}^n \prod_{i = 1}^{k - 1} a^2_{i, j} \right) \left( \sum_{j =
   1}^n a_{k.j}^2 \right) \geqslant \left( \sum_{j = 1}^n \prod_{i = 1}^k
   a_{i, j} \right)^2 . \end{equation}
It follows from (\ref{1}) and (\ref{2}) that
\[ \left( \sum_{v \in E_0} \prod_{i = 1}^k (N_{E_i}^{r_i} (v) - \frac{d}{n} |E_i|)
   \right)^2 \leqslant \prod_{i = 1}^k \sum_{v \in E_0} (|N_{E_i}^{r_i} (v) | -
   \frac{d}{n} |E_i|)^2 \leqslant \lambda^{2 k} \prod_{i=1}^k |E_i|. \]
It can be written as
\begin{equation}\label{3} \left| \sum_{I \subset \{1, \ldots, k\}} (-1)^{k-|I|} \left(\frac{d}{n}\right)^{k -
   |I|} \prod_{j \notin I} |E_j| \sum_{v \in E_0} \prod_{i \in I} N_{E_i}^{r_i} (v) \right| \leqslant
   \lambda^k \sqrt{\prod_{i=1}^{k}|E_i|} . \end{equation}
For any $I \subset \{1, \ldots, k\}$ with $0< |I| < k$, by the induction
hypothesis, we have
\begin{equation}\label{4} \sum_{v \in E_0} \prod_{i \in I} N_{E_i}^{r_i} (v) = e_I (E_0;\{E_i\}_{i \in I}) =(1 + o (1))  \left( \frac{d}{n}
   \right)^{|I|} |E_0|\prod_{i \in I} |E_i|  . \end{equation}
Putting (\ref{3}) and (\ref{4}) together, we have
\[ \left| \sum_{v \in E_0} \prod_{i = 1}^k N_{E_i}^{r_i} (v) - (1 + o (1)) \left( \frac{d}{n}
   \right)^k \prod_{i=0}^k |E_i|  \right| \leqslant \lambda^k \sqrt{\prod_{i=1}^k |E_i|} .
\]
Since $|E_0|^2\prod_{i=1}^k|E_i| \gg ( \frac{n}{d} \lambda)^{2k}$, the left hand side is dominated by $(1 + o (1)) \left( \frac{d}{n} \right)^k \prod_{i=0}^k |E_i| $. This implies that
\[ e_{\{r_1, \ldots, r_k \}} (E_0;\{E1,\ldots,E_k\}) = \sum_{v \in E_0} \prod_{i = 1}^k N_{E_i}^{r_i}
   (v) =  (1 + o (1))\left( \frac{d}{n} \right)^k \prod_{i=0}^k |E_i|, \]
completing the proof of the theorem.  
\end{proof} 

\section{Finite Euclidean and non-Euclidean graphs}
\label{sec:graphs}

In this section, we construct our main tools to study simplexes in vector spaces over finite fields, the graphs associated to finite Euclidean and non-Euclidean spaces. The construction of finite Euclidean graphs follows one of Medrano et al. in \cite{medrano} and the construction of finite non-Euclidean graphs follows one of Bannai, Shimabukuro, and Tanaka in \cite{bannai-shimabukuro-tanaka}. 

\subsection{Finite Euclidean graphs}

Let $\mathbbm{F}_q$ denote the finite field with $q$ elements where $q \gg 1$
is an odd prime power. For a fixed $a \in \mathbbm{F}_q$, the finite Euclidean graph $G_q(a)$ in $\mathbbm{F}_q^d$ is defined as the graph with vertex
set $\mathbbm{F}_q^d$ and the edge set
\[ \{(x, y) \in \mathbbm{F}_q^d \times \mathbbm{F}_q^d \mid x \neq y, ||x - y|| = a\}. \]
In \cite{medrano}, Medrano et al. studied the spectrum of these graphs and showed that these graphs are asymptotically Ramanujan graphs.  Precisely, they proved the following result.

\begin{theorem} \cite{medrano} \label{tool 4}
The finite Euclidean graph $G_q(a)$ is a regular graph with $n(q,a) = q^d$ vertices
of valency

\begin{equation*}
k(q,a) = \left\{
        \begin{array}{ll}
        q^{d-1} + \chi((-1)^{(d-1)/2}a)q^{(d-1)/2} & a \neq 0, \ \ d \ \mbox{odd,}\\
        q^{d-1} - \chi((-1)^{d/2}))q^{(d-2)/2} & a \neq 0, \ \ d \ \mbox{even,}\\
        q^{d-1} & a = 0, \ \ d\ \mbox{odd,}\\
        q^{d-1} - \chi((-1)^{d/2}))(q-1)q^{(d-2)/2} & a = 0,\ \ d \ \mbox{even}. 
    \end{array}  \right.
\end{equation*}
where $\chi$ is the \textit{quadratic character}
\begin{equation*}
\chi(a) = \left\{
        \begin{array}{ll}
        1 & a \neq 0, \ \ a \ \mbox{is square in} \ \mathbbm{F}_q,\\
        -1 & a \neq 0, \ \ a \ \mbox{is nonsquare in} \ \mathbbm{F}_q,\\
        0 & a = 0.
    \end{array}  \right.
\end{equation*}
Let $\lambda$ be any eigenvalues of the graph $G_q(a)$ with $\lambda \neq$ valency of the graph then 
\begin{equation}
|\lambda| \leq 2q^{\frac{d-1}{2}}.
\end{equation}
\end{theorem}

\subsection{Finite non-Euclidean graphs}

Let $V = \mathbbm{F}_q^{2k-1}$ be the $(2k-1)$-dimensional vector space over the finite field $\mathbbm{F}_q$ ($q$ is an odd prime power). For each element $x$ of $V$, we denote the $1$-dimensional subspace containing $x$ by $[x]$. Let $\Omega$ be the set of all square type non-isotropic $1$-dimensional subspaces of $V$ with respect to the quadratic form $Q(x) = x_1^2 + \ldots + x_{2k-1}^2$. The simple orthogonal group $O_{2k-1}(\mathbbm{F}_q)$ acts transtively on $\Omega$, and yields a symmetric association scheme $\Psi(O_{2k-1}(\mathbbm{F}_q),\Omega)$ of class $(q+1)/2$. The relations of $\Psi(O_{2k-1}(\mathbbm{F}_q),\Omega)$ are given by
\begin{eqnarray*}
 R_1 & = & \{([U],[V]) \in \Omega \times \Omega \mid  (U+V) \cdot (U+V) = 0\},\\
 R_i & = & \{([U],[V]) \in \Omega \times \Omega \mid (U+V) \cdot (U+V) = 2 + 2 \nu^{- (i
  - 1)} \} \, (2 \leqslant i \leqslant (q - 1) / 2)\\
 R_{(q+1)/2} & = & \{([U], [V]) \in \Omega \times \Omega \cdot (U+V) \cdot (U+V) = 2\},
\end{eqnarray*}
where $\nu$ is a generator of the field $\mathbbm{F}_q$ and we assume $U\cdot U = 1$ for all $[U] \in \Omega$ (see \cite{bannai-hao-song}).

The graphs $(\Omega,R_i)$ are asymptotic Ramanujan for large $q$. The following theorem summaries the results from Section 2 in \cite{bannai-shimabukuro-tanaka} in a rough form.

\begin{theorem}\cite{bannai-shimabukuro-tanaka}\label{tool5} The graphs $(\Omega,R_i)$ $(1 \leq i \leq (q+1)/2)$ are regular of order $q^{2k-2}(1+o_q(1))/2$ and valency $Kq^{2k-3}$. Let $\lambda$ be any eigenvalue of the graph $(\Omega,R_i)$ with $\lambda \neq$ valency of the graph then 
\[|\lambda| \leq kq^{(2k-3)/2},\]
for some $k,K > 0$ (In fact, we can show that $k = 2+o_q(1)$ and $K = 1+o_q(1)$ or $1/2+o_q(1)$).
\end{theorem}

\section{Proof of Theorem \ref{main2-tes}}
\label{sec:proofs}

We now give a proof of Theorem \ref{main2-tes}. For any $\{a_{ij}\}_{1\leq i<j \leq k+1} \in \mathbbm{F}_q^{\binom{k+1}{2}}$, define
\[
  T_{\{a_{ij}\}_{1\leq i < j \leq k+1}} (E) =\{(x_i)_{i=1}^{k+1} \in E^{k+1} :  \|x_i - x_j\|= a_{ij}\}.
\]

Hart and Iosevich \cite{hart-iosevich} observed that over finite fields, a
(non-degenerate) simplex is defined uniquely (up to translation and rotation)
by the norms of its edges.

\begin{lemma}\label{e-lemma} \cite{hart-iosevich} Let $P$ be a (non-degenerate) simplex with
  vertices $V_0, V_1, \ldots, V_k$ with $V_j \in \mathbbm{F}_q^d$. Let $P'$ be
  another (non-degenerate) simplex with vertices $V_0', \ldots, V_k'$. Suppose
  that
  \begin{equation}
    \|V_i - V_j \|=\|V_i' - V_j' \|
  \end{equation}
  for all $i, j$. Then there exists $\tau \in \mathbbm{F}_q^d$ and $O \in S
  O_d (\mathbbm{F}_q)$ such that $\tau + O (P) = P'$.
\end{lemma}

Therefore, it suffices to show that if $E \subset \mathbbm{F}_q^{2k-1}$, $k \geq 3$, such that $|E| \gg q^{2k-1-\frac{1}{2k}}$, then
\begin{equation}\label{tp2}
  \left|\left\{ \{a_{ij}\}_{1\leq i<j \leq k+1} \in \mathbbm{F}_q^{\binom{k+1}{2}} : |T_{\{a_{ij}\}_{1\leq i < j \leq k+1}} (E) | > 0 \right\} \right| \geqslant cq^{\binom{k+1}{2}}.
\end{equation}

Consider the set of colors $L =\{c_0, \ldots, c_{q - 1} \}$
corresponding to elements of $\mathbbm{F}_q$. We color the complete
graph $G_q$ with the vertex set $\mathbbm{F}_q^{2k-1}$, by $q$ colors such that
$(x, y) \in \mathbbm{F}_q^{2k-1} \times \mathbbm{F}_q^{2k-1}$ is colored by $c_i$ whenever $||x - y|| = i$.  

Suppose that $|E| \gg q^{2k-1-\frac{1}{2k}}$, then we have
\[|E| \gg \left(\frac{q^{2k-1} \cdot 2q^{k-1}}{q^{2k-2}(1+o(1))}\right)^{\frac{2i}{i+2}},\]
for all $2 \leq i \leq k$. From Theorem \ref{tool 4}, $G_q$ is a $(q^{2k-1},q^{2k-2}(1+o(1)), 2q^{k-1})$-r.c. graph when $k \geq 3$. Therefore, applying Theorem \ref{main-tool} for the number of $k$-stars of type $(a_{12},\ldots,a_{1(k+1)})$ in $E^{k+1}$, we have
\[e_{a_{12},\ldots,a_{1(k+1)}}(E;\{E,\ldots,E\}) = \left( \frac{q^{2k-2}(1+o(1))}{q^{2k-1}} \right)^k|E|^{k+1}(1+o(1)) = \frac{|E|^{k+1}(1+o(1))}{q^k},\]
for any $a_{12},\ldots,a_{1(k+1)} \in \mathbbm{F}_q$. 

Let $\mathbbm{F}_{\square}^{*}$ denote the set of non-zero squares in $\mathbbm{F}_q$. For any $a_{12}, \ldots, a_{1(k+1)} \in \mathbbm{F}_{\square}^{*}$, then
\begin{eqnarray*} 
|\{(x_1,\ldots,x_{k+1}) \in E^{k+1}:  \|x_1-x_i\| = a_{1i}\}| & = & e_{a_{12},\ldots,a_{1(k+1)}}(E;\{E,\ldots,E\})\\
& = &\frac{|E|^{k+1}(1+o(1))}{q^k}.\end{eqnarray*}

By the pigeon-hole principle, there exists $x_1 \in E$ such that
\[\{(x_2,\ldots,x_{k+1}) \in E^k: \|x_1-x_i\| = a_{1i}\}| = \frac{|E|^k(1+o(1))}{q^k}.\]

Let $S_t = \{v \in \mathbbm{F}_q^{2k-1} : ||v||=t \}$ denote the sphere of radius $t$ in $\mathbbm{F}_q^{2k-1}$, then $|S_t| = q^{2k-2}(1+o(1))$ for any $t \in \mathbbm{F}_q$. Let $E_i = \{v \in E: \|x_1 - v \| = a_{1i}\} \subset S_{a_i}$, $2 \leq i \leq k+1$, then $|E_2|\ldots|E_{k+1}| = \frac{|E|^k(1+o(1))}{q^k}$ and $|E_2|,\ldots,|E_{k+1}| \leq O(q^{2k-2})$. This implies that
\[
|E_i| \geq \Omega \left( \frac{|E|^k(1+o(1))}{q^{k+(2k-2)(k-1)}} \right) \gg q^{2k-\frac{5}{2}}.
\]

There are $(q-1)^k/2^k$ possible choices for $a_{12}, \ldots, a_{1(k+1)} \in \mathbbm{F}_\square^{*}$, from Lemma \ref{e-lemma}, it suffices to show that $T_k(E_2,\ldots,E_{k+1}) \geq cq^{\binom{k}{2}}$ for some $c>0$. Let $E'_i = \{[x] : x \in E_i\} \subset \Omega$ where $\Omega$ is the set of all square type non-isotropic $1$-dimensional subspaces of $\mathbbm{F}_q^{2k-1}$ with respect to the quadratic form $Q(x) = x_1^2 + \ldots + x_{2k-1}^2$. Since each line through origin in $\mathbbm{F}_q^{2k-1}$ intersects the unit sphere $S_1$ at two points, $|E'_i| \geq |E_i|/2 \gg q^{2k-\frac{5}{2}}$. Suppose that $([U], [V]) \in E'_i \times E'_j$ is an edge of $(\Omega,R_l)$, $2 \leq l \leq (q-1)/2$. Then
\[(U+V)\cdot(U+V) = 2+\alpha_l,\]
where $\alpha_l = 2\nu^{-(l-1)}$. Since $U \cdot U = V \cdot V = 1$, we have $(U-V)\cdot (U-V) = 2 - \alpha_l$. The distance between $U$ and $V$ (in $E'_i \times E'_j$) is either $(U+V)\cdot (U+V)$ or $(U-V) \cdot (U-V)$, so
\begin{equation}\label{tool6} ||U-V|| \in \{2+\alpha_l, 2-\alpha_l\}.\end{equation}

Consider the set of colors $L =\{r_1, \ldots, r_{(q+1)/2} \}$
corresponding to classes of the association scheme $\Psi(O_{2k-1}(\mathbbm{F}_q),\Omega)$. We color the complete graph $P_q$ with the vertex set $\Omega$, by $(q+1)/2$ colors such that
$([U], [V]) \in \Omega \times \Omega$ is colored by $r_i$ whenever $([U],[V]) \in R_i$. 

From Theorem \ref{tool5}, $P_q$ is a $((1+o(1))q^{2k-2}/2,Kq^{2k-3}, kq^{(2k-3)/2})$-r.c. graph when $k \geq 3$. Since $|E'_i| \gg q^{2k-\frac{5}{2}}$, we have
\[
|E'_i| \gg kq^{(2k-3)/2} \left( \frac{(1+o(1))q^{2k-2}/2}{Kq^{2k-3}} \right)^{k-1}.
\]
Therefore, applying Theorem \ref{tool 3} for colored $k$-complete subgraphs
of $P_q$ then $P_q$ contains all possible colored $k$-complete subgraphs. From (\ref{tool6}), $([U],[V])$ is colored by $r_l$ ($2 \leq l \leq (q-1)/2$) then $||U-V|| \in \{2+\alpha_l, 2-\alpha_l\}$. Hence $T_k(E_2,\ldots,E_{k+1}) \geq cq^{\binom{k}{2}}$ for some $c>0$. The theorem follows.

\end{document}